# FINAL-BOUNDARY VALUE PROBLEM IN THE NON-CLASSICAL TREATMENT FOR A SIXTH ORDER PSEUDOPARABOLIC EQUATION

I.G. Mamedov

A.I.Huseynov Institute of Cybernetics of NAS of Azerbaijan. Az 1141, Azerbaijan, Baku st. B. Vahabzade, 9
E-mail: ilgar-mammadov@rambler.ru

**Abstract**

In a rectangular domain we consider a final - boundary value problem for the equation

$$(V_{4,2}u)(x,y) \equiv D_x^4 D_y^2 u(x,y) + a_{3,2}(x,y) D_x^3 D_y^2 u(x,y) +$$
$$+ a_{4,1}(x,y) D_x^4 D_y u(x,y) + \sum_{\substack{i=0 \ j=0 \\ i+j<5}}^{4} \sum_{j=0}^{2} a_{i,j}(x,y) D_x^i D_y^j u(x,y) = Z_{4,2}(x,y).$$

For this equation we consider a final-boundary value problem with non-classical conditions not requiring agreement conditions. Equivalence of these conditions with the classic boundary condition is substantiated in the case if the solution of the stated problem is sought in S.L.Sobolev anisotropic space.

**Keywords:** final boundary value problem, pseudoparabolic equation, discontinuous coefficient equations

## Problem statement

Consider the equation

$$(V_{4,2}u)(x,y) \equiv D_x^4 D_y^2 u(x,y) + a_{3,2}(x,y) D_x^3 D_y^2 u(x,y) +$$
$$+ a_{4,1}(x,y) D_x^4 D_y u(x,y) + \sum_{\substack{i=0 \ j=0 \\ i+j<5}}^{4} \sum_{j=0}^{2} a_{i,j}(x,y) D_x^i D_y^j u(x,y) = Z_{4,2}(x,y) \in L_p(G). \quad (1)$$

Here $u(x,y)$ is a desired function determined on $G$; $a_{i,j}(x,y)$ are the given measurable functions on $G = G_1 \times G_2$, where $G_k = (0, h_k)$, $k = 1,2$; $Z_{4,2}(x,y)$ is a given measurable function on $G$; $D_t^n = \partial^n / \partial t^n$ is a generalized differentiation operator in S.L.Sobolev sense and $D_t^0$ is an identity transformation operator.

Equation (1) is a hyperbolic equation possessing two real characteristics $x = const, y = const$, the first of which is four-fold, the second is two-fold. Therefore, in some sense, we can consider equation (1) as a pseudoparabolic equation [1]. This equation is a generalization of many model equations of some processes (for example, Boussinesq-Liav equation, Manjeron's equation, heat conductivity equation, telegraph equation, string vibration equation and etc.).

In the present paper we consider equation (1) in the general sense when the coefficients $a_{i,j}(x, y)$ are non - smooth functions satisfying only the following conditions:

$$a_{i,j}(x, y) \in L_p(G), \quad i = \overline{0,3} \quad j = \overline{0,1};$$

$$a_{4,j}(x, y) \in L_{\infty,p}^{x,y}(G), \quad j = \overline{0,1}; \quad a_{i,2}(x, y) \in L_{p,\infty}^{x,y}(G), \quad i = \overline{0,3}.$$

Under these conditions the solution $u(x, y)$ of equation (1) will be sought in S.L.Sobolev space

$$W_p^{(4,2)}(G) \equiv \{u(x, y): D_x^i D_y^j u(x, y) \in L_p(G), i = \overline{0,4}, j = \overline{0,2}\},$$

where $1 \leq p \leq \infty$. For equation (1) we can give the final boundary conditions of classic form as (see. Fig.1.)

$$\begin{cases} u(x, y)|_{x=h_1} = \varphi_1(y); & u(x, y)|_{y=h_2} = \psi_1(x); \\ \dfrac{\partial u(x, y)}{\partial x}\bigg|_{x=h_1} = \varphi_2(y); & \dfrac{\partial u(x, y)}{\partial y}\bigg|_{y=h_2} = \psi_2(x); \\ \dfrac{\partial^2 u(x, y)}{\partial x^2}\bigg|_{x=h_1} = \varphi_3(y); & \dfrac{\partial^3 u(x, y)}{\partial x^3}\bigg|_{x=h_1} = \varphi_4(y); \end{cases} \quad (2)$$

where $\varphi_k(y), k = \overline{1,4}$ и $\psi_1(x), \psi_2(x)$ are the given measurable functions on $G$. Obviously, in the case of conditions (2), in addition to the conditions

$$\varphi_k(y) \in W_p^{(2)}(G_2), k = \overline{1,4}; \quad \psi_1(x) \in W_p^{(4)}(G_1), \psi_2(x) \in W_p^{(4)}(G_1);$$

the given functions should satisfy also the following agreement conditions:

$$\begin{cases} \psi_1(h_1)=\varphi_1(h_2); & \psi_2(h_1)=\varphi'_1(h_2); \\ \psi'_1(h_1)=\varphi_2(h_2); & \psi'_2(h_1)=\varphi'_2(h_2); \\ \psi''_1(h_1)=\varphi_3(h_2); & \psi''_2(h_1)=\varphi'_3(h_2); \\ \psi'''_1(h_1)=\varphi_4(h_2); & \psi'''_2(h_1)=\varphi'_4(h_2). \end{cases} \qquad (3)$$

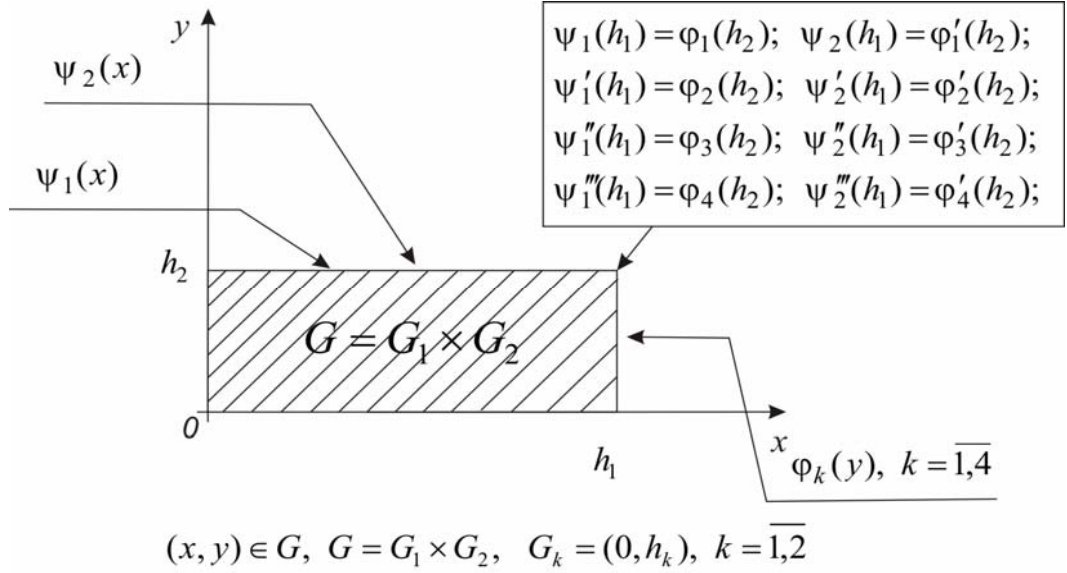

$$(x,y)\in G,\ G=G_1\times G_2,\ G_k=(0,h_k),\ k=\overline{1,2}$$

**Fig.1. Geometrical interpretation of classic final-boundary conditions**

Consider the following non - classical boundary conditions:

$$\begin{cases} V_{i,j}u \equiv D_x^i D_y^j u(h_1,h_2)=Z_{i,j}\in R,\ i=\overline{0,3},\ j=\overline{0,1}; \\ (V_{4,j}u)(x)\equiv D_x^4 D_y^j u(x,h_2)=Z_{4,j}(x)\in L_p(G_1),\ j=\overline{0,1} \\ (V_{i,,2}u)(y)\equiv D_x^i D_y^2 u(h_1,y)=Z_{i,2}(y)\in L_p(G_2),\ i=\overline{0,3} \end{cases} \qquad (4)$$

If the function $u \in W_p^{(4,2)}(G)$ is a solution of classic type final boundary value problem (1), (2), then it is also a solution of problem (1), (4) for $Z_{i,j},\ i=\overline{0,4},\ j=\overline{0,2}$, determined by the following equalities:

$$Z_{0,0}=\varphi_1(h_2)=\psi_1(h_1);\ Z_{0,1}=\varphi'_1(h_2)=\psi_2(h_1);\ Z_{1,0}=\varphi_2(h_2)=\psi'_1(h_1);$$
$$Z_{1,1}=\varphi'_2(h_2)=\psi'_2(h_1);\ Z_{2,0}=\varphi_3(h_2)=\psi''_1(h_1);\ Z_{2,1}=\varphi'_3(h_2)=\psi''_2(h_1);$$
$$Z_{3,0}=\varphi_4(h_2)=\psi'''_1(h_1);\ Z_{3,1}=\varphi'_4(h_2)=\psi'''_2(h_1);\ Z_{4,0}(x)=\psi_1^{(IV)}(x);$$

$$Z_{4,1}(x)=\psi_2^{(IV)}(x); Z_{0,2}(y)=\varphi_1''(y); Z_{1,2}(y)=\varphi_2''(y) Z_{2,2}(y)=\varphi_3''(y);$$

$$Z_{3,2}(y)=\varphi_4''(y).$$

It is easy to prove that the inverse one is also true. In other words, if the function $u \in W_p^{(4,2)}(G)$ is a solution of problem (1), (4), then it is also a solution of problem (1), (2) for the following functions:

$$\varphi_1(y)=Z_{0,0}+(y-h_2)Z_{0,1}+\int_{h_2}^{y}(y-\tau)Z_{0,2}(\tau)d\tau; \tag{5}$$

$$\varphi_2(y)=Z_{1,0}+(y-h_2)Z_{1,1}+\int_{h_2}^{y}(y-\xi)Z_{1,2}(\xi)d\xi; \tag{6}$$

$$\varphi_3(y)=Z_{2,0}+(y-h_2)Z_{2,1}+\int_{h_2}^{y}(y-\eta)Z_{2,2}(\eta)d\eta; \tag{7}$$

$$\varphi_4(y)=Z_{3,0}+(y-h_2)Z_{3,1}+\int_{h_2}^{y}(y-\nu)Z_{3,2}(\nu)d\nu; \tag{8}$$

$$\psi_1(x)=Z_{0,0}+(x-h_1)Z_{1,0}+\frac{(x-h_1)^2}{2!}Z_{2,0}=$$
$$+\frac{(x-h_1)^3}{3!}Z_{3,0}+\int_{h_1}^{x}\frac{(x-\mu)^3}{3!}Z_{4,0}(\mu)d\mu; \tag{9}$$

$$\psi_2(x)=Z_{0,1}+(x-h_1)Z_{1,1}+\frac{(x-h_1)^2}{2!}Z_{2,1}+$$
$$+\frac{(x-h_1)^2}{3!}Z_{3,1}+\int_{h_1}^{x}\frac{(x-\rho)^3}{3!}Z_{4,1}(\rho)d\rho; \tag{10}$$

Note has the functions (5)-(10) possess one important property, more exactly, agreement conditions (3) for all $Z_{i,j}$, having the above-mentioned properties are fulfilled for them automatically. Therefore, we can consider equalities (5) - (10) as a general form of all the functions $\varphi_k(y), k=\overline{1,4}, \psi_1(x), \psi_2(x)$, satisfying agreement conditions (3).

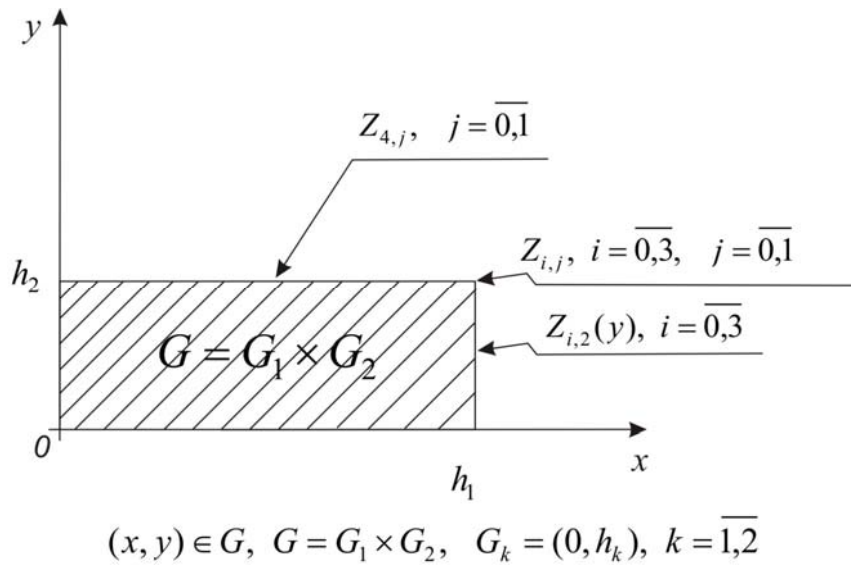

$$(x,y) \in G, \ G = G_1 \times G_2, \ G_k = (0, h_k), \ k = \overline{1,2}$$

**Fig.2. Geometric interpretation of final - boundary conditions in non classical treatment.**

So, the classic form final boundary problem (1), (2) and in non-classical treatment (1), (4) (see fig.2) are equivalent in the general case. However, the final-boundary value problem in non-classical treatment (1), (4) is more natural by the statement than problem (1), (2). This is connected with the fact that in the statement of problem (1), (4) the right sides of boundary conditions don't require additional conditions of agreement type. Note that such boundary value problems in non-classical treatment were considered in the author's papers [2-4].